\documentclass[12pt]{amsart}
\usepackage{latexsym}
\usepackage{graphicx}
\usepackage{amsmath}
\usepackage{amssymb}
\usepackage{verbatim}
\usepackage{amsfonts}

\newtheorem{lemma}{Lemma}
\newtheorem{thm}{Theorem}

\newtheorem*{thma}{Theorem A}
\newtheorem*{thmb}{Theorem B}

\newcommand{\bd}{\partial}

\newcommand{\h}{\mathbb{H}}
\newcommand{\R}{\mathbb{R}}

\begin{document}

\title{Principal curvatures of fibers and Heegaard surfaces}
\author{William Breslin}
\address{Department of Mathematics, University of Michigan, Ann Arbor, MI 48109}
\email{breslin@umich.edu}
\urladdr{http://www-personal.umich.edu/~breslin/index.html}

\begin{abstract}
We study principal curvatures of fibers and Heegaard surfaces smoothly embedded in hyperbolic 3-manifolds. It is well known that a fiber or a Heegaard surface in a hyperbolic 3-manifold cannot have principal curvatures everywhere less than one in absolute value.  We show that given an upper bound on the genus of a minimally embedded fiber or Heegaard surface and a lower bound on the injectivity radius of the hyperbolic 3-manifold, there exists a $\delta > 0$ such that the fiber or Heegaard surface must contain a point at which one of the principal curvatures is greater than $1 + \delta$ in absolute value.
\end{abstract}

\maketitle

\section{Introduction}
The principal curvatures of a surface or lamination smoothly embedded in a hyperbolic 3-manifold are related to the topology of the surface and the 3-manifold. For example in \cite{Breslin} we show that incompressible surfaces and strongly irreducible Heegaard surfaces embedded in hyperbolic 3-manifolds can always be isotoped to a surface with principal curvatures bounded in absolute value by a fixed constant that does not depend on the surface or the 3-manifold.  In \cite{Breslin3} we show that laminations in hyperbolic 3-manifolds with principal curvatures everywhere close to zero have boundary leaves with non-cyclic fundamental group and that laminations in hyperbolic 3-manifolds with principal curvatures everywhere in the interval $(-1,1)$ have boundary leaves with non-trivial fundamental group.

This note was motivated by a question about surfaces with principal curvatures near the interval $(-1,1)$.  It is well known that a closed orientable surface smoothly embedded in a finite-volume complete hyperbolic 3-manifold with principal curvatures everywhere in the interval $(-1,1)$ is incompressible and lifts to a quasi-plane in $\h^3$ (see Thurston's notes \cite{tnotes} or Leininger \cite{Leininger} for a proof).  Thus Heegaard surfaces and fibers in hyperbolic 3-manifolds cannot have principal curvatures everywhere in the interval $(-1,1)$.  We are interested in finding obstructions to isotoping Heegaard surfaces and fibers in hyperbolic 3-manifolds to have principal curvatures close to the interval $(-1,1)$.
See Rubinstein \cite{Rubminimal} or Krasnov-Schlenker \cite{Krasnov} for more on surfaces in hyperbolic 3-manifolds with principal curvatures in the interval $(-1,1)$.



It follows from Freedman-Hass-Scott \cite{FHS} that an incompressible surface in a closed Riemannian 3-manifold can be isotoped to a minimal surface. It follows from work of Pitts-Rubinstein that a strongly irreducible Heegaard surface in a closed Riemannian 3-manifold can be be isotoped to either a minimal surface or the boundary of a regular neighborhood of a minimal surface (see \cite{Rubminimal} for a sketch of the proof).  We show that given an upper bound on the genus of a minimally embedded fiber or Heegaard surface and a lower bound on the injectivity radius of the hyperbolic 3-manifold, there exists a $\delta > 0$ such that the fiber or Heegaard surface must contain a point at which one of the principal curvatures is greater than $1 + \delta$ in absolute value.

\begin{thm}\label{fiber}
For each $g \ge 2$, $\epsilon > 0$, there exists $\delta := \delta(g,\epsilon)$ such that if $S$ is a genus $g$ minimally embedded fiber in a closed hyperbolic mapping torus $M$ with $\operatorname{inj}(M) > \epsilon$, then $S$ contains a point at which one of the principal curvatures is at least $1 + \delta$ in absolute value.
\end{thm}

\begin{thm}\label{heegaard}
For each $g \ge 2$, $\epsilon > 0$, there exists $\delta := \delta(g,\epsilon)$ such that if $S$ is a genus $g$ minimally embedded Heegaard surface in a closed hyperbolic 3-manifold $M$ with $\operatorname{inj}(M) > \epsilon$, then $S$ contains a point at which one of the principal curvatures is at least $1 + \delta$ in absolute value.
\end{thm}

The proofs of Theorem \ref{fiber} and Theorem \ref{heegaard} both use geometric limit arguments.  Assuming that no such $\delta > 0$ exists, we consider a sequence of hyperbolic 3-manifolds as in the statment with minimally embedded fibers or Heegaard surfaces whose principal curvatures are closer and closer to the interval $[-1,1]$.  After possibly passing to a subsequence, the sequence of manifolds converges geometrically to a hyperbolic 3-manifold $M$ and the surfaces converge to an incompressible surface $S$ in $M$ with principal curvatures everywhere in the interval $[-1,1]$.  This implies that the limit set of a lift of $S$ to $\h^3$ is a proper subset of $\bd\h^3$.
In either case, we show that the cover of $M$ corresponding to the image of $\pi_1 (S)$ in $\pi_1 (M)$ has a doubly degenerate hyperbolic structure contradicting that the limit set of a lift of $S$ to $\h^3$ is a proper subset of $\bd\h^3$.

\section{Preliminaries}

Let $M$ be a hyperbolic 3-manifold with no cusps and finitely generated fundamental group.  By a result of Scott, $M$ has a \textit{compact core} which is a compact submanifold $C$ of $M$ whose inclusion into $M$ is a homotopy equivalence.  The connected components of $M \setminus C$ are called the \textit{ends} of $M$.  It follows from the positive solution of the Tameness Conjecture by Agol \cite{agol} and Calegari-Gabai \cite{Cal-Gabai} that an end of $M$ is homeomorphic to $\Sigma \times [0,\infty )$ where $\Sigma$ is a closed orientable surface.
The convex core, $CC(M)$, of $M$ is the smallest convex submanifold of $M$ whose inclusion is a homotopy equivalence.  An end $E$ of $M$ is \textit{convex-cocompact} if $E \cap CC(M)$ is compact and $E$ is \textit{degenerate} otherwise.  Given a closed orientable surface $\Sigma$ of genus greater than one, a hyperbolic structure on $\Sigma \times \R$ such that both ends are degenerate is called \textit{doubly degenerate}.

A sequence of pointed hyperbolic $n$-manifolds $(M_i ,p_i )$ \textit{converges geometrically} to the pointed hyperbolic $n$-manifold $(M,p)$ if for every sufficiently large $R$ and each $\epsilon > 0$, there exists $i_0$ such that for every $i \ge i_0$, there is a $(1 + \epsilon )$-bilipschitz pointed diffeomorphism $\kappa_i : (B(p,R) , p) \rightarrow M_i$, where $B(p,R) \subset M$ is the ball of radius $R$ centered at $p$ and $B(p_i ,R) \subset M_i$ is the ball of radius $R$ centered at $p_i$.  We call the maps $\kappa_i$ \textit{almost isometries}.\\

We will use the fact that minimal surfaces have bounded diameter in the presence of a lower bound on injectivity radius.  See Rubinstein \cite{Rubminimal} or Souto \cite{souto} for more on minimal surfaces in hyperbolic 3-manifolds.

\begin{lemma}\label{mindiam}
Let $S$ be a connected minimal surface in a complete hyperbolic 3-manifold $M$ with $\operatorname{inj}(M) \ge \epsilon$.  Then the diameter of $S$ is at most $4|\chi(F)|/\epsilon + 2\epsilon$.
\end{lemma}

We will also use the following Lemma in the proofs of Theorem \ref{fiber} and Theorem \ref{heegaard}.

\begin{lemma}\label{limitset}
If $S$ is a closed orientable surface smoothly immersed with principal curvatures everywhere in the interval $[-1,1]$ in a complete hyperbolic 3-manifold $M$ with no cusps, then the limit set of a lift of $S$ to $\h^3$ is a proper subset of $\bd\h^3$.
\end{lemma}

\begin{proof}
Let $\tilde{S}$ be a lift of $S$ to $\h^3$.  Assume that $\tilde{S}$ is not a horosphere, as otherwise we are done. Thus the principal curvatures of $S$ cannot be everywhere equal to $1$ or everywhere equal to $-1$. 
If the principal curvatures at every point of $S$ are $-1$ and $1$, then there is a pair of line fields defined on the entire surface, implying that $S$ is a torus.  Since closed surfaces in $M$ with all principal curvatures in $[-1,1]$ are incompressible and $M$ has no cusps, $S$ cannot be a torus.
Thus there is a point $p$ in $\tilde{S}$ at which one of the principal curvatures is in $(-1,1)$.  Assume that the other principal curvature at $p$ is in $[-1,1)$.  Let $H$ be a horosphere tangent to $\tilde{S}$ at $p$. Use an upper half space model of $\h^3$ in which $H$ is a horizontal plane and $\tilde{S}$ is below $H$.  Let $l$ be a simple loop in $\tilde{S}$ which contains $p$ such that the principal curvatures at each point on $l$ are in $[-1,1)$ with at least principal curvature in $(-1,1)$. At each point $x$ in $l$, let $H_x$ be the horosphere above $\tilde{S}$ tangent to $\tilde{S}$ at $x$.  For each $x$ in $l$, let $c_x \in \bd\h^3$ be the center of the horosphere $H_x$.  The set of points $C = \{ c_x | x \in l \}$ forms a closed curve in $\bd\h^3$.  Since the principal curvatures of $\tilde{S}$ are everywhere in the interval $[-1,1]$, $\tilde{S}$ cannot transversely intersect any of the horospheres $H_x$.  Thus, the limit set of $\tilde{S}$ cannot cross the closed curve $C$, so that the limit set of $\tilde{S}$ is a proper subset of $\bd\h^3$.
\end{proof}

It is well-known that the limit set of a lift to $\h^3$ of a fiber $\Sigma$ in a doubly degenerate hyperbolic  $\Sigma \times \R$  is the entire boundary $\bd\h^3$.  By Lemma 2, such a fiber $\Sigma$ cannot be smoothly embedded with principal curvatures everywhere in the interval $[-1,1]$.

\section{Principal curvatures of fibers}

In the proof of Theorem \ref{fiber}, we will use the following well-known fact about geometric limits of hyperbolic mapping tori.

\begin{thma}\label{thma}
Let $(M_i,p_i)$ be a sequence of pairwise distinct pointed hyperbolic mapping tori with genus $g$ fibers and $\operatorname{inj}(M_i) > \epsilon$ for all $i$.  Then a subsequence of $(M_i,p_i)$ converges geometrically to a pointed hyperbolic 3-manifold $(M,p)$ homeomorphic to $\Sigma \times \R$ where $\Sigma$ is a closed genus $g$ surface and $M$ has a doubly degenerate hyperbolic structure.
\end{thma}

\noindent\textit{Proof of Theorem \ref{fiber}.}  Suppose, for contradiction, that Theorem \ref{fiber} does not hold.  Then there exists a sequence of hyperbolic mapping tori $(M_i)$ with $\operatorname{inj}(M_i) > \epsilon$ such that $M_i$ has a genus $g$ minimal surface fiber with principal curvatures less than $1 + 1/i$ in absolute value.  For each $i$, let $p_i$ be a point in $S_i$.  By Theorem A the sequence $(M_i,p_i)$ has a subsequence, say the entire sequence, which converges to a doubly degenerate pointed hyperbolic 3-manifold $(M,p)$ homeomorphic to $\Sigma \times \R$ where $\Sigma$ is a genus $g$ closed surface.  By Lemma \ref{mindiam}, the diameters of the surfaces $S_i$ are uniformly bounded.  Thus we can find a compact subset $K$ of $M$ homeomorphic to $\Sigma \times [-1,1]$ such that for $i$ large enough, say for all $i$, $S_i$ is contained in $\kappa_i (K)$.  
The surface $S := \Sigma\times\{0\}$ in $M$ is isotopic to $\kappa_i^{-1} (S_i)$ for each $i$. 
Since the surfaces $\kappa_i^{-1} (S_i)$ have bounded area and curvature, a subsequence converges to a smoothly immersed surface with principal curvatures in $[-1,1]$ which is homotopic to $S$. Lemma \ref{limitset} implies that the limit set of a lift of $S$ to $\h^3$ is a proper subset of $\bd\h^3$, contradicting the fact that $M$ is doubly degenerate.  \hfill$\Box$



\section{Principal curvatures of Heegaard surfaces}

In the proof of Theorem \ref{heegaard}, we will use the following well-known fact about geometric limits.

\begin{thmb}\label{thmb}
Every sequence $(M_i,p_i)$ of pointed hyperbolic 3-manifolds with $\operatorname{inj}(M_i,p_i)$ bounded away from 0 has a geometrically convergent subsequence.
\end{thmb}

We also need a Lemma from Souto (Lemma 2.1 from \cite{souto2}).

\begin{lemma}\label{sequence}
Let $(M_i)$ be a sequence of hyperbolic 3-manifolds converging to a hyperbolic manifold $M$.  Assume that there is a compact subset $K\subset M$ such that for all sufficiently large $i$ the homomorphism $\pi_1(K) \rightarrow \pi_1(M_i)$ provided by geometric convergence is surjective.  Then, if the cover of $M$ corresponding to the image of $\pi_1(K)$ into $\pi_1(M)$ has a convex-cocompact end, so does $M_i$ for all but finitely many $i$.
\end{lemma}

\noindent\textit{Proof of Theorem \ref{heegaard}.}  Suppose for contradiction that Theorem \ref{heegaard} does not hold.  Then there exists a sequence $(M_i)$ of closed hyperbolic 3-manifolds with $\operatorname{inj}(M_i) > \epsilon$ such that $M_i$ has a genus $g$ minimal Heegaard surface $S_i$ with principal curvatures less than $1 + 1/i$ in absolute value.  For each $i$ let $p_i$ be a point in $S_i$.  By Theorem B the sequence $(M_i,p_i)$ has a convergent subsequence, say the entire sequence, which converges geometrically to a pointed hyperbolic 3-manifold $(M,p)$.  By Lemma \ref{mindiam}, the diameters of the surfaces $S_i$ are uniformly bounded.  Thus each $M_i$ contains a compact subset $K_i$ homeomorphic to $S_i \times [-1,1]$ with uniformly bounded diameter.  For $i$ large enough the pull-back $\kappa_i^{-1}(K_i)$ of $K_i$ through the almost isometries provided by geometric convergence are embedded compact subsets homeomorphic to $\Sigma \times [-1,1]$ where $\Sigma$ is a closed surface of genus $g$.  For $i$ large enough the surfaces $\kappa_i^{-1}(S_i)$ are all isotopic to a fixed embedded genus $g$ surface $S$ in $M$.  
Since the surfaces $\kappa_i^{-1} (S_i)$ have bounded area and curvature, a subsequence converges to a smoothly immersed surface with principal curvatures in $[-1,1]$ which is homotopic to $S$.
Thus the surface $S$ is incompressible in $M$ and by Lemma \ref{limitset} the limit set of a lift of $S$ to $\h^3$ is a proper subset of $\bd\h^3$.

To arrive at a contradiction we will show that the cover of $M$ corresponding to the image of $\pi_1(S)$ into $\pi_1(M)$ is doubly degenerate, implying that the limit set of a lift of $S$ to $\h^3$ is all of $\bd\h^3$.  For $i$ large enough $\kappa_i(S)$ is isotopic to the Heegaard surface $S_i$ in $M_i$, so that the homomorphism $(\kappa_i)_* : \pi_1(S) \rightarrow \pi_1(M_i)$ provided by geometric convergence is surjective.  By Lemma \ref{sequence}, if the cover of $M$ corresponding to the image of $\pi_1(S)$ into $\pi_1(M)$ has a convex-cocompact end, so does $M_i$ for all but finitely many $i$.  Since each $M_i$ is closed we have that the cover of $M$ corresponding to the image of $\pi_1(S)$ into $\pi_1(M)$ cannot have a convex-cocompact end.  Thus the cover of $M$ corresponding to the image of $\pi_1(S)$ into $\pi_1(M)$ is doubly degenerate contradicting the fact that $S$ is isotopic to a surface with principal curvatures everywhere in $[-1,1]$. \hfill$\Box$\\

\textbf{Acknowledgement.}  This work was partially supported by the NSF RTG grant 0602191.

\bibliographystyle{amsalpha}
\bibliography{tri}

\end{document}